\documentstyle[12pt]{article}
\input mssymb.tex
\newtheorem{teo}{Theorem}
\newtheorem{pro}{Proposition}
\newtheorem{definition}{Definition}
\newtheorem{lemma}{Lemma}
\def\C{{\Bbb C}}
\def\Z{{\Bbb Z}}   
\def\R{{\Bbb  R}}
\def\P{{\Bbb  P}}
\def\N{{\Bbb  N}}
\begin{document}
\baselineskip18pt

\bigskip

\bigskip

\bigskip
\bigskip
\title{On the defect of an analytic disc}
\author{ Patrizia Rossi\thanks{Supported by M.P.I. and C.N.R. 
Research Groups.}}
\date{Universita' di Roma - Tor Vergata}
\maketitle
\normalsize
\bigskip

\begin{abstract}
Although the concept of defect of an analytic disc attached to a generic
manifold of $\C^{n}$ seems to play a merely technical role, it turns 
out to be a rather deep and fruitful notion for the extendability of 
CR functions defined on the manifold.

In this paper we give a new geometric description of defect, drawing 
attention to the behaviour of the interior points of the disc 
by infinitesimal perturbations. For hypersurfaces a stronger result
holds because these perturbations describe a complex vector space
of $\C^{n}$.

For a big analytic disc the defect does not need to be smaller than the 
codimension of the manifold. Indeed we show by an example that it can be 
arbitrarily large independently of the codimension of the manifold.\\
Nevertheless we also prove that the defect is always finite.
In the case of a hypersurface we give a geometric upper bound for the defect.
\end{abstract}

\bigskip
\bigskip
\centerline{\bf Introduction.}
\bigskip
\bigskip

The concept of {\sl defect} of an analytic disc attached to a CR manifold
$M\subset\C^{n}$ appeared first in the well known paper of A.E. Tumanov on
the edge of the wedge theorem [T1]. Although the defect seemed there to play
a merely technical role, it turned out later to be a rather deep and fruitful
notion.\\
Besides the study of extendability of CR functions, the defect has interesting
applications to propagation phenomena ([Tr]) and CR-maps ([BR], [CR]).\\
Rather than the original definition, we prefer to start sketching two 
characteristic properties of the defect which are the core of Tumanov's 
theorem.

\medskip
\noindent
Let $M\subset\C^{n}$ be a suitably smooth manifold that we assume to be 
{\sl generic}, i.e. the tangent space $T_{p}M$ to $M$ at any $p\in M$ generates
linearly all of $\C^{n}$ over $\C$.\footnote{ All CR manifolds are locally
CR-equivalent to a generic one.}\\
Let $D$ be the unit disc and $\Gamma$ its boundary. Consider an analytic disc 
$\phi : \bar{D}\rightarrow \C^{n}$ of class $C^{1,\delta}$, $0<\delta<1$,
$\phi\in {\cal O}(D)$, with boundary on $M$ ($\phi\Gamma\subset M$), and
fix a {\sl base point}
$p=\phi (\zeta_{0})$, $\zeta_{0}\in\Gamma$, at the boundary.\\ 
\indent
{\sl Consider first "suitably" small discs.} The infinitesimal perturbations 
$\dot{\phi}$ of $\phi$, keeping fixed the point $\phi(\zeta_{0})$ and still
respecting $\phi\Gamma\subset M$, form a vector space $U$.\\
Now the defect can be described as follows.
\noindent
\begin{itemize}
\item[a)] Fix arbitrarily $\zeta_{1}\in\Gamma$, $\zeta_{1}\neq\zeta_{0}$. As
$\dot{\phi}$ runs through all perturbations in $U$, $\dot{\phi}(\zeta_{1})$
fills a vector subspace of $T_{\phi(\zeta_{1})}M$. The codimension of this 
subspace is the defect of $\phi$.
\noindent
\item[b)] Consider the starting velocity $\vec{v}$ of the curve 
$t\mapsto\phi [(1-t)\zeta_{0}]$ at the base point $p=\phi (\zeta_{0})$
when we move along the radius of the disc
and fix an arbitrary supplementary vector space $S$ to $T_{p}M$ with 
projection $\pi$. As $\phi$ undergoes all perturbations in $U$, 
$\pi\vec{v}$ describes a vector-subspace of $S$ whose codimension is again
the defect of $\phi$. 
\end{itemize}
The main purpose of the present paper is to show that the defect can be also
described by perturbations of the {\sl interior} of the disc in the following
way.
\noindent
\begin{itemize}
\item[c)] Fix arbitrarily an {\sl interior} point $\zeta_{2}\in D$ and 
again subject $\phi$ to the perturbations $\dot{\phi}$. Then this 
differential has an image $V$ whose span over $\C$ has a codimension 
equal to the defect of $\phi$. Furthermore for hypersurfaces 
we obtain a stronger result because $V$ is always a {\sl complex} vector space.
This is the content of our Theorem 2 in section 3.
\end{itemize}
The actual definition of defect in [T1] is neither a) nor b): it is introduced
in a rather technical way (see (\ref{eq:V}) in section 1 below) and apparently
depends on the base point $p$.\\
However this can be easily reformulated , as we did in [R], in order to make 
the defect independent of $p$, and not to require that the disc is small.\\
In section 1 we introduce this "reformulated" definition and prove that 
it is the same as the original one for small discs. This is the same point 
of view of Baouendi- Rothschild-Tr\'epreau in [BRT].\\
\indent
{\sl Let us now speak about large discs.} In this case 
Tumanov's definition of defect has no meaning. We shall use 
our definition and prove that the defect is always
finite.\\
{\sl In any case the characterizations a), b), c) no longer hold for large 
discs.} Indeed, according to those characterizations, the defect 
obviously cannot exceed the codimension of $M$, while in Proposition 1 we show that
the unit disc of the $z_{1}$-axis, viewed as analytic disc attached to a 
particular algebraic real hypersurface, has defect $2k+1$.\\
In Proposition 2 we give an upper bound for the defect of an
analytic disc attached to a hypersurface and in Theorem 1 we 
show that also in higher codimension the defect is finite.\\ 
In order to clarify the geometric construction which leads to 
Theorem 2, which is stated and proved in section 3, we gather in 
section 2 several 
results concerning mainly the Hilbert transform and matrix valued 
functions in the disc. Although they should be considered as a part of the 
proof, some of those results might have some interest in their own.

The author would like to thank E.M. Chirka, A. Huckleberry and N. Krushilin 
for their useful and constructive suggestions.

\bigskip
\bigskip
\centerline{\bf \S 1. The defect of a disc.}
\bigskip
\bigskip

In this section we show, as we did in [R], that the original definition
of defect can be reformulated in a geometric way eliminating the 
particular equation of $M$ and the base point of the disc.\\ 
Let $M$ be a real, {\sl generic} manifold of class $C^{2,\epsilon}$, 
$0<\epsilon <1$ in $\C^{n}$. We shall always assume $M$ to be an open,
relatively compact subset of a larger manifold.\\
The fiber at $p\in M$ of the {\sl holomorphic co-normal bundle} $CM$ of $M$
is the real vector space of the forms $\omega =\sum_{j=1}^{n}a_{j}dz_{j}$
such that ${\cal I}m(\omega)$ vanishes on the tangent space $T_{p}M$.\\ 
If $\{\rho_{1}=\cdots =\rho_{m}=0\}$ are local real equations for $M$, with
$\partial \rho_{1}\wedge\cdots\wedge
\partial\rho_{m}\neq 0$ on $M$ (by the genericity of $M$), then $CM =
i\R\partial\rho_{1}\oplus\cdots\oplus i\R\partial\rho_{m}$.
Let $\phi$ be an analytic disc of class $C^{1,\delta}$, $0<\delta<
\epsilon$, with boundary on $M$, i.e. $\phi\Gamma\subset M$. A 
section of the pull-back $\phi^{*} C\!M\rightarrow \Gamma$ of $CM$ has the 
form $\omega\circ\phi =\sum_{j=1}^{n}(a_{j}\circ\phi)dz_{j}$. We say that
$\omega\circ\phi$ {\sl extends holomorphically into the disc $D$} if all 
coefficients $a_{j}\circ\phi$ extend holomorphically to $D$.
\begin{definition}
The defect $d(\phi)$ of a disc $\phi$ is the dimension of the real vector 
space $$E_{\phi}=\{C^{1,\delta}-sections \ of \ \phi^{*}CM \ which \ extend 
\ holomorphically \ to \ D \}.$$
\end{definition} 
Observe that, if $\alpha$ is an automorphism of the disc, then 
$E_{\phi}=E_{\phi\circ\alpha}$. Thus {\sl the defect is invariant by 
right composition with an automorphism of the disc}.\\
Define the {\sl size} of a disc $\phi$ as 
$$|\phi|=\inf_{z_{0}\in\C^{n}}||\phi +z_{0}||_{1, \delta}.$$
This quantity measures how far is a disc from being a constant disc.\\ 
We shall 
prove in Proposition 3 that, if $|\phi|$ is smaller than a constant
depending only on $M$, then $\phi^{*} C\!M$ {\sl has a moving frame such 
that the sections which extend holomorphically to $D$ have constant 
components with respect to this frame}.
Since $rk(CM)=codim M$, this gives
$$d(\phi)\leq codim M, \ for \ small \ |\phi|.$$
On the other hand Tumanov's characterizations of the defect
(a), b) of the introduction) and our Theorem 2 (or statement c) in the 
introduction) also obviously imply $d(\phi)\leq codim M$.\\ 
We will now show that, for large discs,
$$d(\phi)> codim M$$
can also occur, but in any case the defect is finite.

\medskip

\begin{pro}
The analytic disc $\phi :\zeta\mapsto (\zeta , 0)$, $|\zeta|\leq 1$,
as a disc attached to the real hypersurface
$$M=\{(z_{1},z_{2})\in \C^{2}, \ {\cal R}e(z_{1}^{k}z_{2})=0, \ z_{1}\neq 0\}$$
has defect $2k+1$.
\end{pro}
{\bf Proof.} Setting $\rho(z_{1},z_{2})=2{\cal R}e(z_{1}^{k}z_{2})$, we have
$\rho_{z_{1}}(\zeta ,0)=0$, $\rho_{z_{2}}(\zeta ,0)=\zeta^{k}$. The real 
vector space $E_{\phi}$ of the holomorphic functions in the unit disc,
whose restrictions to $|\zeta|=1$ are equal to $\zeta^{k}$ times a real 
function, obviously coincides with the set of the functions of the type $\zeta^{k}
(p_{k}+\bar{p}_{k})$, where $p_{k}$ is any holomorphic polynomial of degree
$k$. Thus $d(\phi)\equiv dim E_{\phi}=2k+1$. $\Box$

\medskip
\noindent
{\bf Remark.} In the example above, except for its center, the whole analytic 
disc is contained in $M$. This is in fact a concidence. Indeed we only need 
to add the term 
$(|z_{1}|^{2}-1)^{2}$ to the equation of $M$ and have $\phi D\cap M=\emptyset$. 
The functions $\rho_{z_{j}}(\zeta ,0)$, $j=1,2$, as well as the 
defect of $\phi$, will not change.

\medskip
\noindent
\begin{teo}
The defect of an analytic disc attached to a generic manifold is finite.
\end{teo}
{\bf Proof.} If $M$, and hence $\phi^{*}CM\rightarrow \Gamma$ 
have not global equations, we can reduce to this case taking its pull-back 
by the map $\Gamma\rightarrow\Gamma$ defined by $\sigma\mapsto\sigma^{2}$.\\
First we observe that since $M$ is generic the 
complex codimension of $T_{p}^{c}M$ \footnote{$T^{c}_{p}M$ is the complex
tangent space to $M$ at $p$.} equals the real codimension $m$ of $M$.
Now, since $\{T^{c}_{\phi(\sigma)}M, \ \sigma\in\Gamma\}$ is a $C^{1}$ 
family of $m$-codimensional complex vector spaces in $\C^{n}$ which represents
a 0-measure set in the corresponding Grassmannian, indeed there exists 
an open dense set in the Grassmannian where we can choose a $m$-dimensional
complex vector space ${\cal V}$ such that ${\cal V}\cap T^{c}
_{\phi(\sigma)}M=\{0\}$ for all $\sigma\in\Gamma$.\\
After a linear change of coordinates we can assume that ${\cal V}$
is the $\{z_{1},\cdots,z_{m}\}$ plane. Thus the matrix 
$A(\sigma):=(\frac{\partial\rho_{k}[\phi(\sigma)]}
{\partial z_{j}})_{k,j\leq m}$ is non degenerate for all $\sigma\in\Gamma$.
An element of $E_{\phi}$ (see Definition 1) is identified with a 
$\C^{n}$ valued function of
the type $\gamma(\sigma)\rho_{z}[\phi(\sigma)]$, where $\gamma$ and 
$\rho_{z}$ are respectively a real $(1,m)$ and a complex $(m,n)$ matrixes.\\
By definition of $E_{\phi}$, $\gamma\rho_{z}$, and hence $\gamma(\sigma)
A(\sigma)$, extends holomorphically to the disc. Since the matrix $A$ is non
degenerate we shall be done if we prove that if this happens, then 
$b:=\gamma A$, on $\Gamma$, belongs to a finite dimensional vector space.\\
Since $\gamma = bA^{-1}$ is real valued we have 

\begin{equation}
bA^{-1}=\overline{bA^{-1}}, \ on \ \Gamma.
\label{eq:cond}
\end{equation}
Set $u(z):=\overline{b(\frac{1}{\bar{z}})}$. $u$ is a holomorphic 
function for $|z|>1$ and $u|_{\Gamma}=\bar{b}$. Then we can write (\ref{eq:cond}) 
in terms of $u$
\begin{equation}
bA^{-1}=u\overline{A^{-1}}, \ on\ \Gamma.
\label{eq:cond2}
\end{equation}
Now, since $u\in{\cal O}(\C-\bar{D})$ has H$\ddot{o}$lder trace on $\Gamma$,
by Plemelj formula we have
\begin{equation}
\frac{1}{2}u(z)+p.v.\frac{1}{2\pi i}\int_{\Gamma}\frac{u(\zeta)}{\zeta -z}d\zeta
=0, \ \forall z\in\Gamma.
\label{eq:cond3}
\end{equation}
Set $I(u)(z):=p.v.\frac{1}{2\pi i}\int_{\Gamma}\frac{u(\zeta)}{\zeta -z}d\zeta$.
By (\ref{eq:cond2}) we have $u=bA^{-1}\bar{A}$ on $\Gamma$ and substituting
this expression in (\ref{eq:cond3}) and setting $C:=A^{-1}\bar{A}$ we obtain
\begin{equation}
\frac{1}{2}b+I(bC)C^{-1}=0, \ on \ \Gamma.
\label{eq:cond4}
\end{equation}
Analogously, since $b$ is a holomorphic function in $D$, applying again the 
Plemelj formula we have
\begin{equation}
-\frac{1}{2}b+I(b)=0, \ on \ \Gamma.
\label{eq:cond5}
\end{equation}
Subtracting (\ref{eq:cond5}) from (\ref{eq:cond4}) we obtain
that $b$ belongs to the kernel of the Fredholm operator of the second
kind (see [V], p.26) 
$$L(b)(z):=b(z)+p.v.\frac{1}{2\pi i}\int_{\Gamma}b(\zeta)\frac{C(\zeta)C^{-1}
(z)-1}{\zeta -z}d\zeta$$
and thus it varies in a finite dimensional vector space. $\Box$

\medskip
\noindent
In the case of hypersurface one can give a simple geometric upper bound 
for $d(\phi)$ suggested by the example in Proposition 1.

\medskip
\noindent
Consider a complex direction $L$ such that, for all $\sigma\in \Gamma$, 
$T^{c}_{\phi(\sigma)}M\cap L=\{0\}$. The set of such $L$ is open dense in 
$\P^{n-1}(\C)$. Then $l(\sigma)=T_{\phi(\sigma)}\cap L$ is a real straight 
line in $L$ through the origin. When $\sigma$ turns once in the 
circle $\Gamma$, $l(\sigma)$ will turn $k$ times in $L$, $k\in {\Z}$. 
\begin{pro}
Let $k$ be the number of times that the real line $l(\sigma)$ above turns
in the complex direction $L$, while $\sigma$ turns once in the unit circle.
For the defect of $\phi$ we have the upper bound 
$$d(\phi)\leq sup \ (0,-2k+1).$$
\end{pro}

\medskip
\noindent
The proof will be based on the following elementary

\medskip
\noindent
{\bf Lemma} {\sl Let $f:\Gamma\rightarrow \C\setminus\{0\}$ be a $C^{\alpha}$ map,
$0<\alpha<1$, with winding number $s$. Then the real vector space $V_{f}$
of the holomorphic functions $g\in {\cal O}(D)\cap C^{\alpha}(\bar{D})$ such
that $g|_{\Gamma}f^{-1}$ is real, has dimension $sup \ (0,2s+1)$.}

\medskip
\noindent
{\bf Proof.} For $f(\sigma)=\sigma^{s}$ we have obviously $V_{\sigma^{s}}=\{0\}$
when $s<0$, while $V_{\sigma^{s}}$ has $\{\sigma^{s}(\sigma^{j}+
\sigma^{-j}), \ i\sigma^{s}(\sigma^{j}-\sigma^{-j}), \ 0\leq j\leq s\}$
as a basis when $s\geq 0$. \\
Then $dim \ V_{\sigma^{s}}=sup \ (0,2s+1)$. For a general $f$ with winding
number $s$, we can write $f(\sigma)=\sigma^{s}e^{X(s)+iY(s)}$, where $X$, $Y$
are real $C^{\alpha}$ functions on $\Gamma$. Let $T$ be the Hilbert transform.
Set $r=X+TY$ and $h= iY-TY$ so that $h$ is the trace on $\Gamma$ of a 
function in ${\cal O}(D)\cap C^{\alpha}(\bar{D})$. We have 
$f(\sigma)=\sigma^{s}e^{r(\sigma)+h(\sigma)}$ and $g|_{\Gamma}f^{-1}$ is real
if and only if $g|_{\Gamma}e^{-h}\sigma^{-s}$ is real, i.e. $g\in 
e^{h}V_{\sigma^{s}}$. This gives $dim \ V_{f}=dim \ V_{\sigma^{s}}=
sup\ (0,2s+1)$. $\Box$

\medskip
\noindent
{\bf Proof of Proposition 2.} Let $L$ be a complex direction such that
$T^{c}_{\phi(\sigma)}M\cap L=\{0\}$, $\forall\sigma\in \Gamma$. We can assume
that $L$ is the $z_{1}$-axis so that, if $\{\rho =0\}$ is the equation of 
$M$, we have $\rho_{z_{1}}[\phi(\sigma)]\neq 0$, $\forall \sigma\in \Gamma$.
Let $s$ be the winding number of $\sigma\mapsto \rho_{z_{1}}[\phi(\sigma)]$.
(Thus the real line $l(\sigma)$ above turns $k=-s$ times.) The first 
component of a section of $\phi^{*}CM$ is a function of the form 
$r(\sigma)\rho_{z_{1}}[\phi(\sigma)]$ with real $r$. If the section extends
holomorphically to $D$, so will do our function. Thus, by the Lemma, 
$d(\phi)\leq sup\ (0,2s+1)$. $\Box$

\medskip
\noindent
We shall now only deal with small discs and always assume 
\begin{equation}
|\phi|<R, 
\label{eq:R}
\end{equation}
reducing $R$ when it is necessary.

\noindent
In order to refer our definition of defect to the original one, 
it is necessary to give a cartesian form to $M$.\\
First we need the following elementary 
\begin{lemma}
Let $M$ be an open, relatively compact subset of a generic manifold of 
class $C^{2,\epsilon}$ and codimension $m$ in $\C^{n}$.
Then, for every $\lambda>0$, there exists $R(\lambda)>0$ with the following 
properties:\\
for every set $L\subset M$ with $diam L<R(\lambda)$ and $\forall
p\in L$, there exist complex affine coordinates $(z=x+iy,w)\in\C^{m}\times
\C^{n-m}$ with origin at $p$, such that $L$ has a neighbourhood in $M$ which is
contained in the set 
\begin{equation}
x=h(w,y), \ \  \ \  |w|<r, \ \ |y|<r,
\label{eq:M}
\end{equation}
where $h: \{|w|<r, |y|<r\} \longrightarrow I \!\! R^{m}$ 
is a function of class $C^{2,\epsilon}$, satisfying
$h(0,0)=dh(0,0)=0$ and $||h_{y}(w,y)||<\lambda$, $\forall
|w|<r, |y|<r$. Here $||\cdot||$ stands for the matrix norm.
\end{lemma}
The proof is quite standard and so we omit it.

\medskip
\noindent
Now we want to prove that $d(\phi)\leq codim M$, for $R$ sufficiently small 
and show that our definition and Tumanov's original 
definition of the defect are the same. In particular this shows that the 
latter is  independent of the choice of base point.\\ 
For this we need to sketch Tumanov's 
presentation which is rather technical.\\ 
Consider the Hilbert transform $T_{1}: C^{1,\delta}(\Gamma)\rightarrow 
C^{1,\delta}(\Gamma)$ normalized at 1. $T_{1}$ is defined on real functions 
by the fact that $f+iT_{1}f$ extends holomorphically to $D$ and $T_{1}f(1)=0$.
It is a bounded operator and $T_{1}^{2}f=-f$ whenever $f(1)=0$.\\ 
Fix arbitrarily $p\in\phi\Gamma$. Since the defect is invariant by 
automorphisms, we can assume $p=\phi(1)$. Replace $L=\phi\Gamma$ in Lemma 1 
and choose $R<R(1/||T_{1}||)$ in (\ref{eq:R}).\\  
There exists a unique $G:\Gamma\rightarrow GL(m,\R)$
of class $C^{1,\delta}$ such that $\sigma\mapsto G(\sigma)(1+ih_{y}
[\phi(\sigma)])$ extends holomorphically to $D$. In fact, since 
$||h_{y}||<1/||T_{1}||$ on $\phi\Gamma$, we can solve the 
equation 
\begin{equation}
G=1-T_{1}[G(h_{y}\circ\phi)]
\label{eq:G}
\end{equation}
which is equivalent to the holomorphic extendability of $G(1+ih_{y})$ to 
the unit disc and $G(1)=1$.\\
Since $(h_{y}\circ\phi)(1)=0$, we have $T_{1}G=G(h_{y}\circ\phi)$. Again 
using Lemma 1, we can take $R$ in (\ref{eq:R}) so small that
the norm of the matrix $h_{y}\circ\phi$ is 
smaller than an absolute constant
which guarantees that, in addition, the holomorphic extension
of $G(1+i(h_{y}\circ\phi))$ is non-degenerate at all points of $\bar{D}$. 
Indeed, as a solution of a fixed point problem, $G$ depends continuosly
on $h_{y}$ and, for $h_{y}=0$, we have $G=1$, $G(1+i(h_{y}\circ\phi))=1$.\\
The defect of $\phi$, with $|\phi|<R$, was originally defined in [T1] 
as the dimension of the vector space
\begin{equation} 
V_{\phi}=\{c\in \R^{m}| \ \ cG(h_{w}\circ\phi) \ \ extends \ \ holomorphically
\ \ to \ \ D\}.
\label{eq:V}
\end{equation}
The next proposition establishes the identity between our definition of defect
and the original one.
\medskip
\noindent
\begin{pro}
If $\phi\Gamma\subset M$ and $|\phi|$ is smaller than a constant depending on
$M$, then $d(\phi)=dim V_{\phi}$. In particular $d(\phi)\leq m=codim M$
\end{pro}   
{\bf Proof.} As we have seen, a neighbourhood of $\phi\Gamma$ in $M$ is 
contained in the manifold ($\ref{eq:M}$). Thus we can assume that $M$ has 
equation ($\ref{eq:M}$).\\ 
Since $d(\phi)$ does not depend on the equations 
of $M$, we can choose $\rho=x-h(w,y)$ and obtain 
$\partial\rho =\frac{1}{2}[(1+ih_{y})dz -2h_{w}dw]$. In those coordinates 
we have $$E_{\phi}\cong\{\gamma \ | \ \gamma(1+ih_{y}\circ\phi) \ and \ 
 \gamma h_{w}\circ\phi \ \ extend \ \ holomorphically \ \ 
to \ \ D\},$$ 
where $\gamma$ is a real $(m\times m)$-matrix function of class $C^{1,\delta}$
on $\Gamma$.\\
Thus, if $c\in V_{\phi}$, then $\gamma =cG\in E_{\phi}$. If 
viceversa $\gamma\in E_{\phi}$, then $\gamma(1+ih_{y}\circ\phi)$ extends to 
a holomorphic matrix $F$ and $G(1+ih_{y}\circ\phi )$ to a holomorphic,
nondegenerate matrix $g$.\\ 
We have on $\Gamma$ $\gamma G^{-1}=
Fg^{-1}$ and, since the left hand side is real and the right hand side extends 
holomorphically, this is a constant real vector $c$. We obtained $\gamma\in 
E_{\phi} \Rightarrow \gamma=cG$, for some $c\in\R^{m}$. Replacing 
$\gamma =cG$ in the second equation of $E_{\phi}$, we obtain $c\in V_{\phi}$. 
Thus $\gamma\in E_{\phi}\Leftrightarrow \gamma =cG$ with $c\in V_{\phi}$. 
This gives $dim E_{\phi} = dim V_{\phi}$. $\Box$  

\medskip
\noindent
{\bf Remark.} The disc in Proposition 1 can be done arbitrarily small
taking $(\epsilon\zeta ,0)$ instead of $(\zeta ,0)$ without any change in the
conclusion. This seems to be in contradiction with the proposition above.
But the closure of the hypersurface of Proposition 1 is singular, thus this
last cannot be relatively compact in any other manifold, while the manifold 
in Proposition 3 is assumed to have this property.

\bigskip
\bigskip
\centerline{\bf \S 2. Preliminary results.}
\bigskip
\bigskip

We now introduce the Hilbert transform $T_{0}$ normalized at $0$,
i.e. with the condition $\int_{0}^{2\pi}T_{0}f d\theta =0$. Thus 
$T_{0}^{2}f =-f+\frac{1}{2\pi}\int_{0}^{2\pi}f d\theta$. Again using Lemma 1, 
we further reduce $R$ in (\ref{eq:R}) so that $||h_{y}(w,y)||<1/||T_{0}||$.\\
The next lemma gives a relation between the matrix $G$ defined by
(\ref{eq:G}) and the unique solution $G_{0}$ of
\begin{equation}
G_{0}=1-T_{0}[G_{0}(h_{y}\circ \phi)].
\label{eq:G'}
\end{equation}
\begin{lemma}
Let $G$ be the matrix-function defined by {\sf ($\ref{eq:G}$)}. Then there 
exists a constant matrix $C \in GL(m, \R)$ such that $G_{0}=CG$ is the solution
of {\sf ($\ref{eq:G'}$)}. Furthermore $G_{0}(1 +ih_{y})$ extends holomorphically
to $D$ and\\ 
$T_{0}G_{0}=G_{0}(h_{y}\circ\phi) -\frac{1}{2\pi}\int_{0}^{2\pi}G_{0}
(\sigma)h_{y}[\phi(\sigma)]d\theta$, $\sigma =e^{i\theta}$.
\end{lemma}
{\bf Proof.} By the definition of $G$ we have that $G(1+ih_{y})$ and 
$[G(1+ih_{y})]^{-1}$ extend holomorphically to $D$. Then the real matrix 
$C\equiv G_{0}G^{-1}=G_{0}(1+ih_{y})[G(1+ih_{y})]^{-1}$ also extends 
holomorphically to $D$. Thus $C$ is constant.\\
Furthermore, applying to ($\ref{eq:G'}$) the transform $T_{0}$, we obtain 
the last assertion of the lemma. $\Box$

\bigskip
\noindent
Now we give some more results which will be needed in the next section.
\begin{lemma}
Let $f$ be a function of class $C^{\alpha}$, $0<\alpha<1$, with $f(1)=0$. Then 
$$\frac{1}{2\pi} \int_{0}^{2\pi}(f +iT_{1}f)d\theta = -\frac{1}{\pi}
\int_{0}^{2\pi}\frac{f(\sigma)}{\sigma -1}d\theta, \ \ \sigma =e^{i\theta}$$
and $$\frac{1}{2\pi} \int_{0}^{2\pi}(f -iT_{1}f)d\theta = -\frac{1}{\pi}
\int_{0}^{2\pi}\frac{f(\sigma)}{\bar{\sigma} -1}d\theta, 
\ \ \sigma =e^{i\theta}.\footnote
{The integrals converge absolutely because $f(1)=0$.}$$
\end{lemma}
{\bf Proof.} Since both sides of those equalities are linear, it is 
sufficient to prove the lemma 
for real $f$. Observe also that in this case the second part follows 
immediately from the first.\\
Now $f+iT_{1}f$ is the boundary value of a holomorphic function $F$ vanishing
at 1. We must compute $F(0)$. For $\zeta\in\bar{D}$ set $F_{1}(\zeta)=
\frac{1}{2\pi} p.v.\int_{0}^{2\pi}\frac{\sigma +\zeta}{\sigma -\zeta}
f(\sigma)d\theta$, with $\sigma=e^{i\theta}$. On $\Gamma$ we have 
${\cal R}e(F_{1})=f={\cal R}e(F)$. Thus $F(\zeta)=F_{1}(\zeta)-F_{1}(1)=
\frac{\zeta -1}{\pi}\int_{0}^{2\pi}\frac{\sigma f(\sigma)}{(\sigma -\zeta)
(\sigma -1)}d\theta$, $\sigma =e^{i\theta}$. Now, setting $\zeta =0$, we 
have the result. $\Box$

\begin{lemma}
Let $f,g \in C^{\alpha}$, $0<\alpha<1$, be such that $\frac{1}{2\pi}
\int_{0}^{2\pi}f(\sigma)d\theta =0$ and $g(1)=0$. Then 
$$\int_{0}^{2\pi}\frac{fg -(T_{0}f)(T_{1}g)}{\sigma -1} d\theta =0$$
and $$\int_{0}^{2\pi}\frac{fg -(T_{0}f)(T_{1}g)}{\bar{\sigma} -1} d\theta =0.$$
\end{lemma}
{\bf Proof.} As in Lemma 3 it is sufficient to only prove the first equality
for real $f$ and $g$. $f+iT_{0}f$ and $g+iT_{1}g$ are boundary values of 
holomorphic functions $F_{0}$ and $F_{1}$ and, by the hypothesis on $f$,
$F_{0}(0)=0$. On $\Gamma$ we have $F_{0}F_{1}=fg -(T_{0}f)(T_{1}g) +
i(fT_{1}g + gT_{0}f)$ and, since the imaginary part vanishes at 1, we obtain 
$fT_{1}g + gT_{0}f=T_{1}[fg -(T_{0}f)(T_{1}g)]$.
Set $A:=fg -(T_{0}f)(T_{1}g)$ so that $F_{0}F_{1}= A+iT_{1}A$ on $\Gamma$
and $A(1)=0$. Since $F_{0}F_{1}$ vanishes at 0, an application of Lemma 3
to $A$ gives the desired equality. $\Box$

\bigskip
\noindent
\begin{pro}
Let $G_{0}$ be the matrix-function defined by equation {\sf (\ref{eq:G'})},
and let $X$ and $Y$ be functions of class $C^{\alpha}$ on $\Gamma$, 
$0<\alpha<1$, linked by the relation $Y=T_{1}X$. Set 
$K=\frac{1}{2\pi}\int_{0}^{2\pi}
G_{0}(\sigma)h_{y}[\phi(\sigma)]d\theta$, with $\sigma =e^{i\theta}$. Then 
$$\int_{0}^{2\pi}\frac{X-KY}{\sigma -1}d\theta = 
\int_{0}^{2\pi}\frac{G_{0}(X-h_{y}Y)}{\sigma -1}
d\theta$$
and $$\int_{0}^{2\pi}\frac{X-KY}{\bar{\sigma} -1}d\theta = 
\int_{0}^{2\pi}\frac{G_{0}(X-h_{y}Y)}
{\bar{\sigma} -1}d\theta.$$
\end{pro}
{\bf Proof.} It is sufficient to apply the Lemma 4 with $f= G_{0}
-1$ and $g=X$, recalling that $T_{0}G_{0}=G_{0}h_{y}-K$ by Lemma 2. $\Box$

\begin{pro}
Let $g$ be a vector-function of type $1\times s$ depending $C^{1,\delta}$ 
on $\sigma\in\Gamma$.\\
If, for all $f=(f_{1},\ldots,f_{s})\in C^{1,\delta}(\bar{D})\cap {\cal O}(D)$
with $f(0)=f(1)=0$, we have
$$I=\int_{0}^{2\pi}(\frac{a}{\sigma -1}+\frac{\bar{a}}{\bar{\sigma} -1})gf
d\theta =0 \ \ \sigma =e^{i\theta},\footnote{Note that the integrand is
continuous because $gf$ is of class $C^{1,\delta}$ and vanishes at $1$.}$$
with $a\in \C$, then $ag$ and $\bar{a}g$ extend holomorphically into $D$.
\end{pro}
{\bf Proof.} By the arbitrary nature of $f$ we can take $f(\sigma)=
(\sigma -1)\sigma^{l}e_{j}$,
with $l\geq 1$, where $e_{j}$ is a vector of the canonical basis of $\R^{s}$
(note that this $f$ satisfies the conditions $f(0)=0$, $f(1)=0$).\\
Then, since $(\bar{\sigma}-1)^{-1}=-\sigma(\sigma -1)^{-1}$, we obtain that
$(a-\sigma\bar{a})g$ extends into the disc $D$
as a holomorphic function $h$.\\
Suppose $a\neq 0$ and since $a$ is a scalar we obtain that 
$(a-\sigma\bar{a})$ has no zeros into $D$. Then $\frac{h}{a-\sigma\bar{a}}$
is a holomorphic function and so $g$ extends holomorphically. In 
particular $ag$ and $\bar{a}g$ extend holomorphically into $D$. $\Box$

\bigskip
\bigskip
\centerline{\bf \S 3. The main theorem.}
\bigskip
\bigskip

We shall now correctly state and prove the main Theorem. As we saw in Section
1, a disc $\phi$ of class $C^{1,\delta}$ attached to $M$, with $|\phi|<R=R(M)$
has the property that, for each $p\in\phi\Gamma$, $\phi$ is a Bishop-lifting
of a unique analytic disc $w(\zeta)$ lying in the complex tangent space 
$T_{p}^{c}M\equiv \C_{w}^{n-m}$ and the Bishop's lifting maps a neighbourhood
of $w(\zeta)$ in $[{\cal O}(D)\cap C^{1,\delta}(\bar{D})]^{n-m}$ onto a 
neighbourhood of $\phi$ in the set ${\cal M}_{p}$ of the $C^{1,\delta}$ discs 
in $\C^{n}$ satisfying $p\in\phi\Gamma\subset M$, $|\phi|<R$. Since $M$ is 
of class $C^{2,\epsilon}$, with $\delta<\epsilon<1$, the lifting 
$w(\zeta)\mapsto\phi$ is of class $C^{1}$ and thus ${\cal M}_{p}$ has a 
natural structure of a $C^{1}$-manifold. So, refering to the point c) in the
introduction, it makes sense to fix $\zeta_{2}$ in the interior
of $D$ and to differentiate the $C^{1}$ map ${\cal M}_{p}\rightarrow\C^{n}$ 
given by the evaluation at $\zeta_{2}$. Since the group of the automorphisms 
of $D$ acts nicely on the right on ${\cal M}_{p}$ and preserves the defect, we 
can add the condition $\phi(1)=p$ to the discs in ${\cal M}_{p}$.

\bigskip
\noindent
\begin{teo}
The differential of the evaluation map ${\cal M}_{p}\rightarrow
\C^{n}$ given by $\phi\mapsto\phi(\zeta)$ (for fixed $\zeta\in D$), 
has an image $V$ whose span over
$\C$ has complex codimension equal to the defect of the disc $\phi$.\\
For hypersurfaces a stronger result holds because this image is always 
a complex vector space.\\
If and only if the defect is 1, then $\zeta\mapsto V(\zeta)$ as a map
$D\rightarrow \P^{n-1}(\C^{n})$ is a holomorphic extension of the map
$\Gamma\rightarrow\P^{n-1}(\C^{n})$ given by $\sigma\mapsto 
T^{c}_{\phi(\sigma)}M$.\\
Furthermore we have $d(\phi)=0$ if and only if $V(\zeta)=\C^{n}$ for one 
(and thus all) $\zeta\in D$.
\end{teo}
{\bf Proof.} By previous discussion we can consider as point 
of evaluation the point $\zeta =0$.\\
Taking coordinates at $p$ as in Lemma 1, with the restrictions 
we imposed in the statement of Proposition 3, the element of ${\cal M}_{p}$
corresponding to $w(\zeta)$ will be $(w(\zeta),z(\zeta))$, where $z(\zeta)
=x(\zeta)+iy(\zeta)$ is defined by its boundary value $z(\sigma)$, $\sigma
\in\Gamma$, with $x(\sigma)$ determined by
\begin{equation}
x(\sigma)=h(w(\sigma),y(\sigma)), \ \ \sigma=e^{i\theta},
\label{eq:wz}
\end{equation}
and $y(\sigma)$ is uniquely determined by the 
Bishop's equation
\begin{equation}
y(\sigma)=T_{1}x(\sigma)=T_{1}h(w(\sigma),y(\sigma)), \ \ \sigma=e^{i\theta}.
\label{eq:y}
\end{equation}

\bigskip
\noindent
From Poisson's formula we have 
$$z(\zeta)=\frac{1}{2\pi }\int_{\Gamma}\frac{\sigma +\zeta}{\sigma -\zeta}
x(\sigma)d\theta +iy(0), \ \ \sigma =e^{i\theta}, \ \ |\zeta|<1.$$
When $|\zeta|=1$, the integral must be taken in the sense of a principal value.\\
We must differentiate the composed map
$$w(\zeta)\mapsto (w(\zeta),z(\zeta))\mapsto (w(0),z(0)),$$
which is defined by ($\ref{eq:wz}$) and ($\ref{eq:y}$) and where $w(\zeta)$
will vary in the Banach space $W$ of vector functions 
$w=w(\zeta): \bar{D}\rightarrow \C^{s}$ of class $C^{1,\alpha}(\bar{D})$, holomorphic
in $D$, with the property $w(1)=0$.\\ 
From Lemma 3 and ($\ref{eq:wz}$) we obtain
$$z(0)=-\frac{1}{\pi}\int_{0}^{2\pi}
\frac{x(\sigma)}{\sigma -1}d\theta=-\frac{1}{\pi}\int_{0}^{2\pi}
\frac{h(w(\sigma),y(\sigma))}{\sigma -1}d\theta, \ \ \sigma =e^{i\theta}.$$
Now we differentiate this expression with respect to $w(\zeta)\in W$, taking
($\ref{eq:wz}$) and ($\ref{eq:y}$) into account.\\
If dot means the differentiation with respect to $w(\zeta)$, on 
$\Gamma$ we have 
$$\left\{\begin{array}{l}
\dot{x}=h_{w}\dot{w} + h_{\bar{w}}\dot{\bar{w}} + h_{y}\dot{y}\\
\dot{y} =T_{1}\dot{x}
\end{array}
\right.$$
where $\dot{x}$, $\dot{y}$ depend $\R$-linearly on $\dot{w}$. If we set
$X$ and $Y$ for their $\C$-linear parts, we have $\dot{x}=X+\bar{X}$, 
$\dot{y}=Y+\bar{Y}$ with
$$\left\{\begin{array}{l}
X=h_{w}\dot{w} + h_{y}Y\\
Y=T_{1}X.
\end{array}
\right.$$
Consider now the case $m=1$.\\
For $a\in\C$, $b\in\C^{n-1}$ we set 
$$l(a,b,\dot{w})=\C-linear \ part \ of \ a\dot{z}(0)+\bar{a}\dot{\bar{z}}
(0)+ b\dot{w}(0)+\bar{b}\dot{\bar{w}}(0)$$
where "$\C$-linear" refers to the dependence on $\dot{w}\in W$.\\
A real subspace is a complex subspace of complex codimension $d$ if and only 
if its annihilator is complex subspace of dimension $d$.\\
So we must only prove that the space 
$$A\equiv\{(a,b)\in\C\times\C^{n-1} \ \ | \ \ l(a,b,\dot{w})=0, \ \
\forall\dot{w}\in W\}$$
is a complex space of dimension $d(\phi)$.\\
The linear parts of the restrictions to $\Gamma$ of $\dot{z}$ and
$\dot{\bar{z}}$ are $X+iY$ and $X-iY$. Thus we have 
$$l(a,b,\dot{w})=\frac{1}{2\pi}\int_{0}^{2\pi}[a(X+iY) + \bar{a}(X-iY)]d\theta
+b\dot{w}(0).$$
Setting $a'=a(1+iK)^{-1}$, where $K$ is defined in Proposition 4, we have\\
$[a(X+iY) + \bar{a}(X-iY)]=a'[(X-KY) + iT_{1}(X-KY)] + \bar{a}'
[(X-KY) - iT_{1}(X-KY)]$, because $T_{1}Y=-X$. If we apply now 
the Lemma 3, we obtain
$$l(a,b,\dot{w})=-\frac{1}{\pi} \int_{0}^{2\pi}
(\frac{a'}{\sigma -1} + \frac{\bar{a}'}{\bar{\sigma} -1})(X-KY)
d\theta +b\dot{w}(0).$$
Thus, using the Proposition 4 for $X=h_{w}\dot{w}+h_{y}Y$ and 
setting $a''=Ca'$, where $C$ is defined in Lemma 2, we have
$$l(a,b,\dot{w})=-\frac{1}{\pi} \int_{0}^{2\pi}
(\frac{a''}{\sigma -1} + \frac{\bar{a}''}{\bar{\sigma} -1})G\ h_{w}\dot{w}
d\theta +b\dot{w}(0).$$
Now we choose $\dot{w}$ such that $\dot{w}(0)=0$ and apply
the Proposition 5 with $g=Gh_{w}$ and $f=\dot{w}$, and obtain that 
$l(a,b,\dot{w})$ vanishes for all such $\dot{w}$ if and only if $a''\in 
V_{\phi}+iV_{\phi}$, where $V_{\phi}$ is defined by (\ref{eq:V})
in section 1. Hence $a\in C^{-1}(V_{\phi}+iV_{\phi})(1+iK)$ and 
this is a $\C$-linear space L having complex dimension $dim_{\R}V_{\phi}= 
d(\phi)$. Thus, if $(a,b)\in A$, then $a\in L$.\\
We now assume $a\in L$ and attempt to find $b$.\\
Since $\frac{1}{\bar{\sigma}-1}=-\frac{\sigma}{\sigma -1}$, the equation 
$l(a,b,\dot{w})=0$ can be written in the form 
\begin{equation}
b\dot{w}(0) = \frac{1}{\pi}\int_{0}^{2\pi}\frac{a''-\sigma\bar{a}''}{\sigma -1}
G\ h_{w}\dot{w}d\theta, \ \ \sigma =e^{i\theta}. 
\label{eq:b}
\end{equation} 
We choose $\dot{w}=(1-\sigma)e_{j}\in W$, where $e_{j}$ is the canonical
basis of $\R^{n-1}$.\\
Substituting this expression in ($\ref{eq:b}$), we obtain 
$$b=\frac{i}{\pi}\int_{\Gamma}(\frac{a''}{\sigma}-
\bar{a}'')G\ h_{w}d\sigma.$$
On the other hand, since $(\bar{a}''+a'')G\ h_{w}$ extends holomorphically
to $D$, we have 
$\int_{\Gamma}\bar{a}''G\ h_{w}d\sigma=
-\int_{\Gamma}a''G\ h_{w}d\sigma$. Then $b$ is given by
$$b = \frac{a''i}{\pi}\int_{\Gamma}\frac{\sigma +1}{\sigma}G\  
h_{w}d\sigma = \frac{i}{\pi}Ca(1+iK)^{-1}\int_{\Gamma}\frac{\sigma +1}{\sigma}
G \ h_{w}d\sigma$$
and this is a $\C$-linear function of $a\in L$. Thus $A$, as the graph of 
a $\C$-linear function on the complex $d(\phi)$-dimensional space $L$, is 
itself a complex $d(\phi)$-dimensional space.

\medskip
\noindent
In the general case ($m>1$) we take $a\in\C^{m}$, $b\in\C^{n-m}$ and set
$$\lambda_{1}(a,b,\dot{w})=\C-linear \ part \ of \ a\dot{z}(0)+ b\dot{w}(0),$$
$$\lambda_{2}(a,b,\dot{w})=\C-linear \ part \ of \ \bar{a}\dot{\bar{z}}(0)+
\bar{b}\dot{\bar{w}}(0)$$
where "$\C$ - linear" refers again to the dependence on $\dot{w}\in W$.\\
As above we have to prove that the space 
$$B\equiv\{(a,b)\in\C^{m}\times\C^{n-m} \ \ | \ \ \lambda_{j}(a,b,\dot{w})=0, 
\ \ \forall\dot{w}\in W \ and \ j=1,2\}$$
has dimension $d(\phi)$.\\
The expressions of $\lambda_{j}(a,b,\dot{w})$ are given by
$$\lambda_{1}(a,b,\dot{w})=\frac{1}{2\pi}\int_{0}^{2\pi}a(X+iY)d\theta
+b\dot{w}(0),$$
$$\lambda_{2}(a,b,\dot{w})=\frac{1}{2\pi}\int_{0}^{2\pi}\bar{a}(X-iY)d\theta.$$
Repeating the previous computations we obtain
$$\lambda_{1}(a,b,\dot{w})=-\frac{1}{\pi} \int_{0}^{2\pi}
\frac{a''}{\sigma -1}G\ h_{w}\dot{w}d\theta +b\dot{w}(0),$$
$$\lambda_{2}(a,b,\dot{w})=-\frac{1}{\pi} \int_{0}^{2\pi}
\frac{\bar{a}''}{\bar{\sigma} -1}G\ h_{w}\dot{w}d\theta.$$
If we choose $\dot{w}$ such that $\dot{w}(0)=0$, we have $a''\in V_{\phi}
+iV_{\phi}$ and hence also in this case $a\in L$.\\
Now to show the dependence on $a$ of $b$ it is enough to choose in 
$$b\dot{w}(0) = \frac{1}{\pi}\int_{0}^{2\pi}\frac{a''}{\sigma -1}
G\ h_{w}\dot{w}d\theta, \ \ \sigma =e^{i\theta}$$
$\dot{w}=(1-\sigma)e_{j}$ obtaing, as above, that $b$ is a 
$\C$-linear function of $a\in L$. $\Box$

\bigskip

\centerline{\bf \S 4. A counterexample in the case of higher codimension.}

\bigskip
\noindent
In this section we will show that the Theorem 2 is not true if 
the codimension of the manifold is greater than $1$.\\
Consider a manifold $M\subset\C^{3}$ of real codimension $2$, having
equation 
$$x=h(w)$$
with $h\in C^{\infty}(\C,\R^{2})$, $h(0)=dh(0)=0$.\\
Call ${\cal V}_{\phi}$ the image of the differential of the function
${\cal M}_{p}\ni\phi\mapsto \phi(0)\in\C^{3}$, $p\in M$
\footnote{Recall that ${\cal M}_{p}$ is the $C^{1}$ Banach manifold of suitably 
small analytic discs attached to $M$ through $p$.}.\\
We will show the existence of a sequence of analytic discs 
$\phi_{\nu}\in {\cal M}_{p}$ with $d(\phi_{\nu})=0$, 
$||\phi_{\nu}||_{1,\delta}\rightarrow 0$ and ${\cal V}_{\phi_{\nu}}
\neq\C^{3}$.

\begin{pro}
There exists a function $h$ defining the manifold $M$ with the 
properties described above and 
$h(\frac{\sigma -1}{\nu})=0$, for $|\sigma|=1$ and $ \nu\in\N$, such 
that\\
(i) the disc $\phi_{\nu}(\zeta)=(0,\frac{\zeta -1}{\nu})\in {\cal M}_{p}$ 
has defect 0;\\
(ii) ${\cal V}_{\phi_{\nu}}\neq\C^{3}$.
\label{pro:esh}
\end{pro}
For the proof of the proposition we need the following

\begin{lemma}
Set $\Gamma_{\nu}=\{\frac{\sigma -1}{\nu}\in\C, \ |\sigma|=1\}$
and let $r_{\nu}\in C^{\infty}(\Gamma_{\nu},\R^{2})$ be given for 
each $\nu\in\N$.\\
Then there exists $h$ as above and for each 
$\nu$, $f_{\nu}\in C^{\infty}(\Gamma_{\nu},\R)$, $f_{\nu}\not\equiv 0$ 
such that we have 
$$h=0, \  \ D_{\rho}h=f_{\nu}r_{\nu}, \ on \ \Gamma_{\nu}$$
where $D_{\rho}$ stays for the normal derivatives to 
$\Gamma_{\nu}$.
\label{lemma:lesh}
\end{lemma}
{\bf Proof.} For fixed $\nu$ we choose $F_{\nu}\in C^{\infty}
(\R^{2})$, $0\leq F_{\nu}\leq 1$, which vanishes only on
$\bigcup_{j\neq\nu}\Gamma_{j}\cup\{0\}$ and $F_{\nu}\equiv 1$ 
out of a very big compact. If necessary, $F_{\nu}$ can be 
changed in a small neighborhood of the point $-2/\nu$ in order
to make sure that $D_{\rho}F_{\nu}\not\equiv 0$ on  
$\Gamma_{\nu}$. Fix also $v_{\nu}\in C^{\infty}(\C,\R^{2})$
bounded and such that
$$v_{\nu}=0, \ D_{\rho}v_{\nu}=r_{\nu}, \ on \ \Gamma_{\nu}.$$
If we choose a real sequence $\lambda_{\nu}$ very rapidly 
decreasing to $0$, then 
$$h=\sum_{\nu=1}^{\infty}\lambda_{\nu}F_{\nu}
v_{\nu}$$ converges obviously to a smooth function which will be 
our function. Indeed among the $F_{j}$'s, 
only $F_{\nu}$ is non vanishing on $\Gamma_{\nu}$, but $v_{\nu}$ 
vanishes there, thus $h|_{\Gamma_{\nu}}=0$.\\
So we only need to take $f_{\nu}=\lambda_{\nu}
D_{\rho}F_{\nu}$. Also $h$ vanishes with its gradient at $0$ because 
so does each $F_{\nu}$ (indeed $F_{\nu}(0)=0$ and 
$F_{\nu}\geq 0$). $\Box$

\bigskip
\noindent
{\bf Proof of the Proposition \ref{pro:esh}.} Set
$a=(1,i)\in\C^{2}$ and choose in the
Lemma \ref{lemma:lesh} $r_{\nu}(\sigma)=(|1+\sigma|^{2}
,2{\cal I}m \sigma)=(1+\bar{\sigma})(a+\sigma\bar{a})$, $|\sigma|=1$.\\
For $\sigma=e^{i\theta}$ we have\\
$0=D_{\theta}h(\frac{\sigma -1}{\nu})=\frac{i}{\nu}\sigma h_{w}
(\frac{\sigma -1}{\nu})-\frac{i}{\nu}\bar{\sigma} h_{\bar{w}}
(\frac{\sigma -1}{\nu})$\\
thus\\
$D_{\rho}h(\frac{\sigma -1}{\nu})=\frac{1}{\nu}\sigma h_{w}
(\frac{\sigma -1}{\nu})+\frac{1}{\nu}\bar{\sigma} h_{\bar{w}}
(\frac{\sigma -1}{\nu})=\frac{2}{\nu}\sigma h_{w}
(\frac{\sigma -1}{\nu})$.\\
Hence on $\Gamma_{\nu}$ 
\begin{equation}
h_{w}(\frac{\sigma -1}{\nu})=\frac{\nu}{2}\bar{\sigma}D_{\rho}h=
\frac{\nu(1+\bar{\sigma})}{2}f_{\nu}(\sigma)(\bar{\sigma}a+\bar{a}).
\label{eq:contr}
\end{equation}
For proving {\sl (i)} we assume that
$ch_{w}(\frac{\sigma -1}{\nu})$ extends holomorphically into $D$ for
some $c=(c_{1},c_{2})\in\R^{2}$.
Set $C=c \ ^{t}a=c_{1}+ic_{2}$. We obtain that the scalar function
$$f_{\nu}(\sigma)(1+\bar{\sigma})(C\bar{\sigma}+\bar{C})=\frac{f_{\nu}
(\sigma)}
{\sigma^{2}}(1+\sigma)(C+\bar{C}\sigma)$$
extends holomorphically.\\
If $C$ is not zero then $C+\bar{C}\sigma$ only vanishes on 
$\Gamma$ and thus $\frac{f_{\nu}(\sigma)}{\sigma^{2}}$ 
extends holomorphically. But this is impossible because 
$f_{\nu}$ is real and not zero.\\
Thus $C$, and consequently $c$, vanishes. This proves {\sl (i)}.\\
For proving {\sl (ii)} it is sufficient to prove that the form 
$\omega=adz+\bar{a}d\bar{z}\not\equiv 0$ vanishes on 
${\cal V}_{\phi_{\nu}}$, $\forall \nu$.\\
We have $<\omega, (\dot{z}(0),\dot{w}(0))>=a\dot{z}(0)+
\bar{a}\dot{\bar{z}}(0)$.\\
Now, for $\phi(\zeta)=(z(\zeta),w(\zeta))$, we have
$$z(0)=-\frac{1}{\pi}\int_{0}^{2\pi}\frac{h[w(\sigma)]}{\sigma -1}d\theta,$$
and thus
$$\dot{z}(0)=-\frac{1}{\pi}\int_{0}^{2\pi}\frac{h_{w}\dot{w}+h_{\bar{w}}
\dot{\bar{w}}}{\sigma -1}d\theta, \ \sigma=e^{i\theta}.$$
Therefore we can write\\
$a\dot{z}(0)+\bar{a}\dot{\bar{z}}(0)=-\frac{1}{\pi}\int_{0}^{2\pi}
\frac{(a-\sigma\bar{a})h_{w}[w(\sigma)]}{\sigma -1}\dot{w}d\theta-
\frac{1}{\pi}\int_{0}^{2\pi}
\frac{(a-\sigma\bar{a})h_{\bar{w}}[w(\sigma)]}{\sigma -1}\dot{\bar{w}}
d\theta$.\\
Passing to $\phi_{\nu}$, we take $w(\sigma)=\frac{\sigma-1}{\nu}$ and 
by (\ref{eq:contr})
obtain $(a-\sigma\bar{a})h_{w}[w(\sigma)]=0$ and 
$(a-\sigma\bar{a})h_{\bar{w}}[w(\sigma)]=0$ because
$(a-\sigma\bar{a})(\bar{a}+\bar{\sigma}a)=0$. Thus 
$a\dot{z}(0)+\bar{a}\dot{\bar{z}}(0)=0$. $\Box$

\begin{quote}
\raggedright{\sl Dip. Matematica}\\
\raggedright{\sl Univ. Tor Vergata}\\
\raggedright{\sl Via della Ricerca Scientifica}\\
\raggedright{\sl 00133 Roma - Italy}
\end{quote}


\bigskip

 


\begin{thebibliography}{Courant}


\bibitem[BR]{BR} M.S. Baouendi - L.P. Rothschild, {\sl A generalized
complex Hopf lemma and its applications to CR mappings.}
Invent. Math. {\bf 111} (1993), 331-348.



\bibitem[BRT]{BRT} M.S. Baouendi - L.P. Rothschild - J.M. Tr\'epreau, {\sl 
On the geometry of analytic discs attached to real manifolds.}
J. Diff. Geom. (to appear).



\bibitem[CH]{CH} R. Courant - D. Hilbert, {\sl Methods of mathematical 
physics.} vol.2 Interscience Publisher, New York (1962).

  

\bibitem[CR]{CR} E.M.Chirka - C. Rea, {\sl Normal and tangent ranks of 
CR mappings.} Duke Math. J. Vol. 76, N.2 (1994), 417-431.



\bibitem[R]{R} P. Rossi, {\sl Estensione in un wedge di funzioni CR.}
Tesi di laurea. Univ. di Roma Tor Vergata (1991), 53 pp.


\bibitem[T1]{T1} A.E. Tumanov, {\sl Extension of CR functions into a wedge
from a manifold of finite type.} Math. Sbornik {\bf 178} (1988), 128-139.
English Tansl. in Math. USSR Sbornik {\bf 64} (1989), 129-140.

\bibitem[T2]{T2} A.E. Tumanov, {\sl Extension of CR functions into a wedge.}
Math. Sbornik {\bf 181} (1990), 385-398. English Transl. in Math. USSR 
Sbornik {\bf 70} (1991), 385-398. 


\bibitem[Tr]{Tr} J.M. Tr\'epreau, {\sl Sur la propagation des singularit\'es
dans les vari\'et\'es CR.} Bull. Soc. Math. Fr. {\bf 118} (1990), 403-450.


\bibitem[V]{V} N.P. Vekua, {\sl Systems of singular integral equations.}
P. Noordhoff Ltd. Groningen, The Netherlands (1967).

\end{thebibliography}
\end{document}